\documentclass[11pt,leqno]{article}
\usepackage{graphicx, amsfonts, amsthm, amsxtra, verbatim, multicol}
\usepackage[mathscr]{euscript}
\begin{document}

\begin{center}
{\LARGE\bf MPHS I: Differential modular forms,Elliptic curves and
Ramanujan foliation}
\\
\vspace{.25in} {\large {\sc Hossein Movasati}} \\
Technische Universit\"at Darmstadt \\
Fachbereich Mathematik \\
Schlo\ss gartenstr. 7,
64289 Darmstadt, Germany \\
Email: {\tt movasati@mathematik.tu-darmstadt.de}
\end{center}

\newtheorem{theo}{Theorem}
\newtheorem{exam}{Example}
\newtheorem{coro}{Corollary}
\newtheorem{defi}{Definition}
\newtheorem{prob}{Problem}
\newtheorem{lemm}{Lemma}
\newtheorem{prop}{Proposition}
\newtheorem{rem}{Remark}
\newtheorem{conj}{Conjecture}
\newcommand\diff[1]{\frac{d #1}{dz}} 
\def\End{{\rm End}}              
\def\hol{{\rm Hol}}
\def\sing{{\rm Sing}}            
\def\spec{{\rm Spec}}            
\def\cha{{\rm char}}             
\def\Gal{{\rm Gal}}              
\def\jacob{{\rm jacob}}          
\newcommand\Pn[1]{\mathbb{P}^{#1}}   
\def\Z{\mathbb{Z}}                   
\def\Q{\mathbb{Q}}                   
\def\C{\mathbb{C}}                   
\def\as{\mathbb{U}}                  
\def\ring{{R}}                         
\def\R{\mathbb{R}}                   
\def\N{\mathbb{N}}                   
\def\A{\mathbb{A}}                   
\def\uhp{{\mathbb H}}                
\newcommand\ep[1]{e^{\frac{2\pi i}{#1}}}
\newcommand\HH[2]{H^{#2}(#1)}        
\def\Mat{{\rm Mat}}              
\newcommand{\mat}[4]{
     \begin{pmatrix}
            #1 & #2 \\
            #3 & #4
       \end{pmatrix}
    }                                
\newcommand{\matt}[2]{
     \begin{pmatrix}                 
            #1   \\
            #2
       \end{pmatrix}
    }
\def\ker{{\rm ker}}              
\def\cl{{\rm cl}}                
\def\dR{{\rm dR}}                

\def\hc{{\mathsf H}}                 
\def\Hb{{\cal H}}                    
\def\GL{{\rm GL}}                
\def\pedo{{\cal P}}                  
\def\PP{\tilde{\cal P}}              
\def\cm {{\cal C}}                   
\def\K{{\mathbb K}}                  
\def\k{{\mathsf k}}                  
\def\F{{\cal F}}                     
\def\M{{\cal M}}
\def\RR{{\cal R}}
\newcommand\Hi[1]{\mathbb{P}^{#1}_\infty}
\def\pt{\mathbb{C}[t]}               
\def\W{{\cal W}}                     
\def\Af{{\cal A}}                    
\def\gr{{\rm Gr}}                
\def\Im{{\rm Im}}                
\newcommand\SL[2]{{\rm SL}(#1, #2)}    
\newcommand\PSL[2]{{\rm PSL}(#1, #2)}  
\def\Res{{\rm Res}}              

\def\L{{\cal L}}                     
\def\Aut{{\rm Aut}}              
\def\any{R}                          
\newcommand\ovl[1]{\overline{#1}}    

\def\pm{{\mathsf p}{\mathsf m}}      
\def\T{{\cal T }}                    
\def\tr{{\mathsf t}}                 
\newcommand\mf[2]{{M}^{#1}_{#2}}     
\newcommand\bn[2]{\binom{#1}{#2}}    
\def\ja{{\rm j}}                 
\def\Sc{\mathsf{S}}                  
\newcommand\es[1]{g_{#1}}            
\newcommand\V{{\mathsf V}}           
\newcommand\Ss{{\cal O}}             
\def\rank{{\rm rank}}                

\def\Ra{\mathrm{Ra}}
\def\nf{g_2}                         
\begin{abstract}
In this article we define the algebra of differential modular forms and we
prove that it is generated by Eisenstein series of weight
$2,4$ and $6$. We define Hecke operators on them, find some analytic relations
between these Eisenstein series and obtain them in a natural way as
coefficients of a family of elliptic curves. Then we
describe the relation between the dynamics of a foliation in $\C^3$ induced
by the Ramanujan relations,  with  vanishing of elliptic integrals.
The fact that a complex manifold over the Moduli of Polarized
Hodge Structures in the case $h^{10}=h^{01}=1$ has an algebraic structure with
an action of an algebraic group plays a basic role in all of the proofs.
\end{abstract}
{\tiny \tableofcontents }
\section{Introduction}
In \cite{gr70} Griffiths introduced the Moduli of Polarized Hodge
Structures/the period domain $D$ and described
a dream to enlarge $D$ to a moduli space of degenerating polarized
Hodge structures. Since in general $D$ is not a Hermitian symmetric domain,
he asked for the existence of a certain automorphic cohomology theory for
$D$, generalizing the usual notion of automorphic forms on symmetric
Hermitian domains. Since then there have been many efforts in the first
part of Griffiths's dream (see \cite{kaus00} and references there) 
but the second part still lives in darkness.

I was looking for some analytic spaces over $D$ for which one may state
Baily-Borel theorem on the unique algebraic structure of quotients
of symmetric Hermitian domains by discrete arithmetic groups
(see \cite{mil03, babo66}).  I realized that even
in the simplest case of Hodge structures, namely  $h^{01}=h^{10}=1$,  
such spaces are
not well studied. This led me to the definition of a new class of
holomorphic functions 
on the Poincar\'e upper half plane which generalize the classical modular
forms. Since a differential operator acts on them we call them
differential modular forms. These new functions are no longer interpreted
as holomorphic sections of a positive line bundle on some compactified moduli
curve. Nevertheless, they appear in a natural way as coefficients in 
families of elliptic curves, analogus to  Eisenstein series in 
the Weierstrass Uniformization Theorem. 
Due to this, I realized that there are very natural
holomorphic foliations  in the coefficient space of families of varieties
whose dynamics is in close relation with the
abelian integrals (resp. Hodge structures) of the family.
In the case of a three parametric elliptic curve, the mentioned foliation is
called Ramanujan foliation because he was the first who obtained a particular
leaf of this foliation using Eisenstein series
(see \cite{la95} Chapter X and (\ref{raman}) bellow).
This article will touch some aspects of number theory,
holomorphic foliations and Hodge theory which
I am going to explain below:

{\bf Number theory:}
Recall the Eisenstein series
\begin{equation}
\label{eisenstein}
\es{k}(z)=a_k{\Big (}1+(-1)^k\frac{4k}{B_k}\sum_{n\geq
1}\sigma_{2k-1}(n)e^{2\pi i z n}{\Big )},\ \  k=1,2,3, \ z\in\uhp
\end{equation}
where $B_k$ is the $k$-th Bernoulli number ($B_1=\frac{1}{6},\
B_2=\frac{1}{30},\ B_3=\frac{1}{42},\ \ldots$), $\sigma_i(n):=
\sum_{d\mid n}d^i$,
\begin{equation}
\label{akhar}
a_1=2\zeta(2)\frac{-1}{2\pi i},\
a_2=2\zeta(4)\frac{60}{(2\pi i)^2},\  a_3=2\zeta(6)\frac{-140}{(2\pi i)^3}
\end{equation}
and $\uhp:=\{x+iy\in \C \mid y>0\}$ is the Poincar\'e upper half plane.
The most well-known differential modular form, which is not a differential
of a  modular form, is the Eisenstein series $\es 1$.
The idea of differentiating modular forms and getting new modular forms
is old and goes back to Ramanujan.
Nesterenko's method (see \cite{nes01}) for the proof of transcendency 
properties of certain numbers is also based on differential
equations satisfied by modular forms.
However, the precise definition of differential modular forms has been given
recently in \cite{bu95}. In the present article we give another 
slightly different definition of differential modular forms (see \S \ref{defi}) over a modular subgroup $\Gamma \subset \SL 2\Z$.
It is based on a canonical behavior of holomorphic functions on the
Poincar\'e upper half plane under the action of $\SL 2\Z$. This
approach has the advantage that it can be generalized to any modular
subgroup of $\SL 2\Z$ but the one in \cite{bu95} works only in the case of
full modular group $\SL 2\Z$.
The set of differential modular forms in the present article is a 
bigraded algebra
$\mf{}{}=\oplus_{n\in\N_0,m\in \N}\mf{n}{m}$, $\mf{0}{m}$ being the set
of classical modular forms of weight $m$, in which the differential
operator $\diff{}$ maps $\mf{n}{m}$ to $\mf{n+1}{m+2}$. We have 
$\es{1}\in \mf{1}{2}, \es{2}\in \mf{0}{4},\es{3}\in\mf{0}{6}$ and we prove:
\begin{theo}
\label{23feb05}
 The functions $g_1,g_2,g_3$ are algebraically
independent and  $\mf{}{}$ is freely generated by
$\es{1},\es{2}$ and $\es{3}$ as a $\C$-algebra.
In particular each $\mf nm$ is a finite
dimensional $\C$-vector space.
\end{theo}
This theorem generalizes the first theorem in each modular forms book that
the algebra of modular forms is generated freely by the Eisenstein series $g_2$
and $g_3$.  Our proof  gives us also the Ramanujan
relations  between the $\es{i}$'s. We define the  action of
Hecke operators on $\mf{n}{m}$ and it turns out that this is similar to the
case of modular forms:
\begin{equation}
\label{quedia}
T_pf(z)=p^{m-n-1} \sum_{d\mid p, 0\leq b\leq d-1}   
d^{-m}f{\Big (}\frac{nz+bd}{d^2}{\Big )},\ p\in \N,\ f\in \mf{n}{m}.
\end{equation}
Hecke operators of this type appear in particular  in the study of
the transfer operator from statistical mechanics which plays an important
role in the theory of dynamical zeta functions (see \cite{HMM}).
The differential operator commutes with Hecke operators (\S \ref{algebra})
and so it induces a map from the set of
new differential forms (see \S \ref{oldnew}) to itself.
Another result which we prove in this article and which 
could be interesting from the number 
theory point of view is the following: Let
$$
g:=(g_1,g_2,g_3): \uhp\rightarrow \C^3
$$
and
$$
T:=\C^3\backslash \{(t_1,t_2,t_3)\in\C^3 \mid 27t_3^2-t_2^3=0\}.
$$
\begin{theo}
\label{realanal}
There are unique analytic functions
$$
B_1, B_2:\ T\rightarrow \R,\  B_3:T\rightarrow \C
$$
such that $B_1$ does not depend on the variable $t_1$ and
\begin{equation}
\label{B1}
B_1\circ g(z)=\Im(z),\ B_1(t_1,t_2k^{-4},t_3k^{-6})=B_1(t)|k|^2
\end{equation}
\begin{equation}
\label{B2}
B_2\circ g=0,\ B_2(
 t_1k^{-2}+k'k^{-1},
t_2k^{-4}, t_3k^{-6})=B_1(t)|k'|^2+B_2(t)|k^{-1}|^2+
\Im(B_3(t)k'\overline{k^{-1}})
\end{equation}
\begin{equation}
\label{B3}
B_3\circ g=1,\ B_3 (
 t_1k^{-2}+k'k^{-1},
t_2k^{-4}, t_3k^{-6})=
B_3(t)k\overline{k^{-1}}+2\sqrt{-1}k
\overline{k'}B_1(t)
\end{equation}
 for all $k\in\C^*$ and $k'\in \C$. Moreover, $|B_3|$ restricted to
 the zero locus of $B_2$ is identically one.
\end{theo}

{\bf Holomorphic singular foliations:} Nowadays the theory of holomorphic
singular foliations is getting a part of Algebraic Geometry.
For a general background in this theory see
(\cite{casa87, camo, bru00,  mov0, CLS89}).
In this direction this article has  two novelties which together exhibit a
connection between the area of holomorphic foliations and
Arithmetic Algebraic
Geometry.
The first one is as follows:
After calculating the Gauss-Manin connection of the following family of 
elliptic curves
\begin{equation}
\label{khodaya}
{\cal E}_t: y^2-4t_0(x-t_1)^3+t_2(x-t_1)+t_3,\ t\in\C^4
\end{equation}
and  considering its relation with the inverse of the period map, we get the
following ordinary differential equation:
\begin{equation}
\label{raman}
\Ra: \left \{ \begin{array}{l}
\dot t_1=t_1^2-\frac{1}{12}t_2 \\
\dot t_2=4t_1t_2-6t_3 \\
\dot t_3=6t_1t_3-\frac{1}{3}t_2^2
\end{array} \right.
\end{equation}
which is called the Ramanujan relations/differential equation/foliation,
because he had observed that $g$ is a solution of (\ref{raman})
(one gets  the classical relations by changing the coordinates
$(t_1,t_2,t_3)\mapsto (\frac{1}{12}t_1,\frac{1}{12}t_2,
\frac{2}{3(12)^2}t_3)$). We denote by
$\F(\Ra)$ the singular holomorphic foliation induced by (\ref{raman}) in
$\C^3$. Its singularities
$$
\sing(\Ra)=\{(t_1,12t_1^2,8t_1^3)\mid t_1\in \C \}
$$
form  a one-dimensional curve in $\C^3$. Consider the family (\ref{khodaya}) with
$t_0=1$ and define
$$
K:={\Big \{}t\in T\mid \int_{\delta}\frac{xdx}{y}=0,\ \hbox{ for some }
\delta\in H_1({\cal E}_t,\Z){\Big \}}
$$
and
$$
M_r:=\{t\in T \mid B_2(t)=r\},\ r\in\R. 
$$
The last part of Theorem \ref{realanal} says that $|B_3|$ restricted
to $M_0$ is identically $1$. We also define 
$$
N_w:=\{t\in M_0\mid B_3(t)=w \}, \ |w|=1, \ w\in\C.
$$
For $t\in \C^3\backslash\sing(\F(\Ra))$ we denote by $L_t$ the leaf
of $\F(\Ra)$ through $t$.
\begin{theo}
\label{foli}
The following is true:
\begin{enumerate}
\item
 The real analytic varieties $M_r,r\in\R,\ N_w, |w|=1$, and the
set $K$ are $\F(\Ra)$-invariant.
\item
The set $K$ is a dense subset of $M_0$ with the following property:
For all $t\in K$ the leaf $L_t$ intersects $\sing(\Ra)$
transversally at some point $p$. 
\item
For all $t\in T$ the leaf $L_t$ has an accumulation point at $T$ if
and only if $t\in M_0$.
\end{enumerate}
\end{theo}
The item $2$ says that there is a transverse disk to $\sing (\F(\Ra))$ at
some point $p$ such that $D\backslash\{p\}$ is a part of the leaf $L_t$.  
For the proof of the above theorem and a precise description of the
foliation $\F(\Ra)$ see \S \ref{ramanfoli}.
The proof of this theorem is based on the fact that the foliation
$\F(\Ra)$ restricted to $T$ is  uniformized by the inverse of the 
period map (see for instance \cite{li00} for similar topics). 

The second novelty in the area of holomorphic foliations is as follows:
\begin{theo}
\label{ellip}
There is no elliptic curve $E$ and a differential form of the second type
$\omega$ on $E$, both defined over $\overline{\Q}$, such that
\begin{equation}
\label{abel}
\int_\delta\omega=0
\end{equation}
for some non-zero topological cycle $\delta\in
H_1(E,\Z)$.
\end{theo}
This theorem  uses Nesterenko's Theorem (see \cite{nes01}) on transcendence
properties of the values of Eisenstein series.
The above theorem for the case in which $\omega$ is of the first kind, 
is well-known. In this case we can even state it for the field $\C$. 
However, it is trivially false when $\omega$ is a differential form of the 
second kind  and we allow transcendental coefficients in $\omega$ or the 
elliptic curve. Generalizations of such theorems to arbitrary curves can be 
derived from the Abelian Subvariety Theorem (see \cite{bo94}, \cite{hovi}).
Abelian integrals of the type (\ref{abel}) appear in deformations of
holomorphic foliations with a first integral in complex surfaces 
(see \cite{ga02, mov0}). The above result shows that the zeros 
of arithmetic abelian integrals are quite different from  complex ones. 

{\bf Hodge theory:}  This article stimulates the hope to realize the second
part of Griffiths' dream with a different formulation.
Differential modular forms can also be introduced for
a complex manifold $\pedo$ over Griffiths  period domain $D$ with an action 
of an algebraic group $G_0$
from the right. Since they are no longer interpreted 
as sections of positive line
bundles over moduli spaces, the question of the existence of a kind of
Baily-Borel
Theorem for $\pedo$ arises. 
In the case of Hodge structures with $h^{01}=h^{10}=1$ we have
\begin{equation}
\label{perdomain}
\pedo:= \left \{
\mat {x_1}{x_2}{x_3}{x_4}\in \GL(2,\C)\mid \Im(x_1\ovl{x_3})>0 \right \},\ 
D=\uhp 
\end{equation}
\begin{equation}
\label{alggroup}
G_0=\left \{\mat{k_1}{k_3}{0}{k_2}\mid \ k_3\in\C, k_1,k_2\in \C^*\right \}
\end{equation}
and we show in this article that $\SL 2\Z\backslash \pedo$ has a canonical
structure of an algebraic quasi-affine variety such that the action of $G_0$
from the right is algebraic.  More precisely, we prove that
$\SL 2\Z\backslash \pedo$ is biholomorphic to
$\C^4\backslash \{t=(t_0,t_1,t_2,t_3)\in \C^4\mid 27t_0t_3^2-t_2^3=0\}$ and under this
biholomorphism the action of $G_0$ is given by:
\begin{equation}
\label{action}
t\bullet g:=(t_0k_1^{-1}k_2^{-1},
 t_1k_1^{-1}k_2+k_3k_1^{-1},
t_2k_1^{-3}k_2, t_3k_1^{-4}k_2^2) $$ $$
 t=(t_0,t_1,t_2,t_3)\in\C^4,
g=\mat {k_1}{k_3}{0}{k_2}\in G_0
\end{equation}
The mentioned biholomorphism is given by the period map (see \S \ref{permap}).

{\bf Singularity theory:}
The Brieskorn module/lattice and its Gauss-Manin connection in singularity
theory  play the role of de Rham cohomology of fibered varieties,
and it is a useful
object when one wants to calculate the Gauss-Manin connection.
In \S \ref{singular} we introduce the associated Brieskorn
$\C[t,\frac{1}{t_0}]$-module $H$ of the family (\ref{khodaya}).
We prove that $H$ is freely generated and we calculate the
Gauss-Manin connection in a canonical basis of  $H$.
We generalize a classical Weierstrass Theorem and see that for the inverse
of the period map, $t_i$ appears as the Eisenstein series $\es{i},\ i=1,2,3$.
The novelty is the appearance of $t_1$ as the Eisenstein series of weight $2$.
This makes us think about a generalization of K. Saito's primitive form
theory (see \cite{sai01}), to the case in which the deformation of a
singularity is bigger than the versal deformation.

{\bf Final note:} An elliptic curve (beside $\Pn 1$) is one of the most
well-studied objects in (Arithmetic) Algebraic Geometry.
This makes me feel that some parts of this work has connections to the works
of many authors that I have not mentioned here. Here I express my sincere
apologizes. My aim of writing this article was just to show, by some simple
examples, how the dynamics of foliations and algebraic geometry of fibered
varieties can be connected.
One of the points in this article which I enjoyed very much, 
was the calculation of the Gauss-Manin connection and the Ramanujan relations 
by {\sc Singular}, \cite{GPS01}. 
The corresponding algorithms were developed in the
articles \cite{mo, hos005} and the corresponding library in {\sc Singular} is
called {\tt brho.lib}. Using this library for  families of algebraic varieties one can obtain certain holomorphic foliations, whose
dynamics has to do with the Hodge structure of Algebraic varieties and in
particular with abelian integrals.
This opens a new connection between (Arithmetic) Algebraic Geometry and the
theory of holomorphic foliations.

Let us now explain the structure of this article. 
\S \ref{dmf} is devoted to
the definition of differential modular forms and the action of Hecke operators
on them. This section is independent of the other sections.
We shall only use the property (\ref{g2}) of the Eisenstein series
$\es{1}$. The reader who is interested only in \S \ref{dmf}, may lose
our proof for the main theorem of this section (Theorem \ref{23feb05}).
\S \ref{singular} is devoted to calculation of the Gauss-Manin
connection of the family (\ref{khodaya}). This section is based on
the machinery introduced in \cite{mo}.
The heart of the present paper is \S \ref{permap} in which we
prove that the period map $\pm$ is a biholomorphism and then we take
the inverse of $d\pm$ to obtaine the Ramanujan relations.
In \S \ref{ramanfoli} we have described the dynamics of the Ramanujan 
foliation and in particular its uniformization by the inverse of the period
map in some quasi-affine variety.
Finally, \S \ref{proofs} is devoted to 
the proofs announced in the Introduction.


{\bf Acknowledgement:}
The main ideas of this paper took place in my mind when
I was visiting Prof. Sampei Usui
at Osaka University. Here I would like to thank him for encouraging me to
study Hodge theory and for his help to understand it.
I would like to thank Prof. Karl-Hermann Neeb for  his interest and 
carefull reading of the present article.

\section{$\mf nm$-functions}
\label{dmf}
In this section we use the notations
$A=\mat {a_A}{b_A}{c_A}{d_A}\in\SL 2\R$ and
$$
T=\mat 1101,\ Q=\mat 0{-1}10,\  x=\mat {x_1}{x_2}{x_3}{x_4},\ 
g=\mat {k_1}{k_2}{0}{k_3},\  x,g\in \GL(2,\C).
$$
When there is no confusion we will simply write  $A=\mat abcd$.
We denote by $\uhp$ the Poincar\'e
upper half plane and
$$\ja (A,z):=c_Az+d_A$$ For $A\in \SL 2\R$ and $m\in\Z$ we use the
slash operator
$$
f|_mA=(\det A)^{m-1}\ja(A,z)^{-m}f(Az)
$$
For a ring $R$ we denote by $\Mat_p(2,R)$ the set of $2\times 2$-matrices
in $R$ with the determinant $p$.
\subsection{Definitions}
\label{defi}
In this section we define the notion of an $\mf nm$-function. For $n=0$
an $\mf 0m$-function is a classical modular form of weight
$m$ on $\uhp$ (see bellow). For arbitrary $n$ we define it by
induction: A holomorphic function $f$ on $\uhp$ is called $\mf nm$
if the following three conditions are satisfied:
\begin{enumerate}
\item
\begin{equation}
\label{4feb05}
f|_m A-f=\sum_{i=1}^{n}\bn{n}{i} c_A^i\ja(A,z)^{-i}f_i, \ \forall A\in \SL 2\Z 
\end{equation}
where $f_i$ is an $\mf {n-i}{m-2i}$-function.
\item
By induction we can write
\begin{equation}
\label{21feb05}
 f_i|_{m-2i} A-f_i=\sum_{j=1}^{n-i}\bn{n-i}{j}
c_A^j\ja(A,z)^{-j}f_{ij}
\end{equation}
where $f_{ij}$ is an $\mf {n-i-j}{m-2i-2j}$-function.
We assume that  $f_{ij}=f_{i-1,j+1}=\ldots=f_{1,i+j-1}=f_{i+j}$.
\item
$f$ has finite growth when $\Im(z)$ tends to $+\infty$, i.e.
$$
\lim_{\Im(z)\to +\infty}f(z)=a_\infty <\infty,\ a_\infty\in\C
$$
\end{enumerate}
The above definition can be made using a subgroup $\Gamma\subset \SL 2\Z$. In 
this article we mainly deal with full differential modular forms, 
i.e. the  case $\Gamma=\SL 2\Z$. 
We will also denote  by $\mf nm$ the set of $\mf nm$-functions and
we set
$$
\mf{}{}:=\oplus_{n,m\in\Z, n\geq 0}\mf nm
$$
A classical modular form satisfies 1 and 2, where instead of the
right hand side of (\ref{4feb05}) we write $0$. Note that for an
$\mf nm$-function $f$ the associated functions $f_i$ are unique
(consider the right hand side of (\ref{4feb05}) as a polynomial in
$c_A\ja (A,z)^{-1}$ with coefficients $\bn{n}{i}f_{i}$). The
second condition guarantees that the consequent use of
(\ref{4feb05}) for two matrices $A,B\in \SL 2\Z$ leads to the same
result. To see this, it is useful to define
\begin{equation}
\label{26feb05} f||_mA:=(\det A)^{m-n-1}\sum_{i=0}^{n}\bn{n}{i}
c_{A^{-1}}^i\ja(A,z)^{i-m}f_{i}(Az),\ A\in \GL (2, \R),\ f\in\mf nm
\end{equation}
The factor $\det A$ is introduced because of Hecke operators
(see \S \ref{heop}).
This is not an action of $\GL (2,\R)$ on $\mf nm$ from the right. The
equalities (\ref{4feb05}) and (\ref{21feb05}) are written in the form
\begin{equation}
\label{13feb05} f_i=f_i||_{m-2i}A:=
\end{equation}
$$
(\det A)^{m-n-i-1}\sum_{j=0}^{n-i}\bn{n-i}{j}
c_{A^{-1}}^j\ja(A,z)^{j+2i-m}f_{i+j}(Az),
\  i=0,1,2,\ldots, n, \ f_0:=f
$$
for all $A\in \SL 2\Z$ (We have substituted $A^{-1}z$ for $z$ and then
$A^{-1}$ for $A$).
\begin{lemm}
\label{toulouse}
We have
$$
f||_m A=f||_m(BA),\ \forall A\in \GL(2,\R),\ B\in\SL 2\Z
$$
\end{lemm}
\begin{proof}
The term $(\det A)^{n-m+1}f||_mA (z)$ is equal to:
\begin{eqnarray*}
&= & \sum_{i=0}^{n}\bn{n}{i} c_{A^{-1}}^i\ja(A,z)^{i-m}f_i(B^{-1}BAz) \\
& = &
\sum_{i=0}^{n}\sum_{j=0}^{n-i}\bn{n}{i}\bn{n-i}{j}
c_{A^{-1}}^ic_{B^{-1}}^j\ja(A,z)^{i-m}\ja(B^{-1},BAz)^{m-2i-j}f_{i+j}(BAz) \\
&= &
\sum_{r=0}^{n}\sum_{s=0}^{r}
\bn{n}{s}\bn{n-s}{r-s}
c_{A^{-1}}^sc_{B^{-1}}^{r-s}\ja(A,z)^{s-m}\ja(B^{-1},BAz)^{m-r-s}f_{r}(BAz) \\
&= & \sum_{r=0}^{n} \bn{n}{r} \ja(BAz, z)^{r-m}f_r(BAz)
\ja(A,z)^{-r}(\sum_{s=0}^{r}\bn{r}{s}\ja(BA,z)^{s}c_{A^{-1}}^s c_{B^{-1}}^{r-s}) \\
&=& \sum_{r=0}^{n} \bn{n}{r} \ja(BA, z)^{r-m}f_r(BAz)
\ja(A,z)^{-r}(\ja(BA,z)c_{A^{-1}}+c_{B^{-1}})^{r} \\
&=& \sum_{r=0}^{n} \bn{n}{r} \ja(BA,
z)^{r-m}c_{(BA)^{-1}}^rf_r(BAz)= (\det A)^{n-m+1}(f||_mBA)(z)
\end{eqnarray*}
In the second line we have used (\ref{4feb05}). 
In the third line we have have changed the counting parameters:
$i+j=r,\ i=s,\ 0\leq s\leq r$. In the other lines we have used
$$
\ja(AB,z)=\ja(A,Bz)\ja(B,z)
$$
and
$$
\ja(BA,z)c_A+\det(A)c_B=c_{BA}\ja(A,z),\  \forall A,B\in\GL(2,\R)
$$
\end{proof}
Since $f|_mT=f$ we can write the Fourier expansion of $f$ at
infinity
$$
f=\sum_{n=-N}^{+\infty}a_n q^{n},\ a_n\in \C,\ N=0,1,2,\ldots,\infty,\
\ q=e^{2\pi i z}.
$$
The third condition on $f$ implies that $N=0$.

\subsection{Algebra of $\mf nm$-functions}
\label{algebra}
In this section we use for $A\in\SL 2\Z$ the simplification $c=c_A$.
Recall the Eisenstein series (\ref{eisenstein}) and
\begin{equation}
\label{jdelta}
\Delta(z):=(27g_3^2(z)-g_2^3(z))=-(\frac{2\pi i}{12})^3q\prod_{n=1}^{\infty}
(1-q^n)^{24}=q-24q^2+252q^3+\cdots
\end{equation}
$$
j(z):=\frac{g_2^3(z)}{-\Delta(z)}=q^{-1}+744+
196884 q+\cdots
$$
Note that $\zeta(2)=\frac{\pi^2}{6},\zeta(4)=\frac{\pi^4}{90},\zeta(6)=
\frac{\pi^6}{945}$ and so
\begin{equation}
\label{pinfty}
p_\infty:=(a_1,a_2,a_3)=(\frac{2\pi i}{12},12(\frac{2\pi i}{12})^2 ,
8(\frac{2\pi i}{12})^3)
\end{equation}
where $a_i$'s are defined in (\ref{akhar}). 
For $k\geq 2$ one can write
$$
\es {k}(z)=s_k\sum_{0\neq (m,n)\in\Z^2}\frac{1}{(n+mz)^k}\in \mf 0k
$$
where $s_2=\frac{60}{(2\pi i)^2}$ and $s_3=\frac{-140}{(2\pi i)^3}$.
Now $\es{1}$ satisfies
\begin{equation}
\label{g2}
\es{1}\mid_2A-\es{1}=c\ja(A,z)^{-1},\ A\in\SL 2\Z
\end{equation}
and so $\es{1}\in \mf 12$ (see for instance \cite{apo90} p. 69).
The following proposition describes the algebraic structure 
of $\mf {n}{m}$:
\begin{prop}
\label{sarita}
The followings are true:
\begin{enumerate}
\item
For an $f\in \mf 1m$ the function $z(f(\frac{-1}{z})-f(z))$
is in $\mf 0{m-2}$, i.e. it is a modular form of weight $m-2$.
\item
$\mf 12$ is a one dimensional $\C$-vector space generated by $g_1$.
\item
If $n\leq n'$ then $\mf nm\subset \mf {n'}{m'}$ and
$$
\mf nm\mf {n'}{m'}\subset \mf {n+n'}{m+m'}, \ \mf nm+\mf {n'}{m}=\mf
{n'}m
$$
\item
For a modular form $f$ of weight $m$ we have $f(g_1)^n\in \mf
n{2n+m}$.
\end{enumerate}
\end{prop}
\begin{proof}
The first item is a direct consequence of the definition applied to $A=Q$.
The only modular forms of weight $0$ are constant functions. This and 1
imply that $\mf 12$ is one dimensional.   Item 3 is derived from the
definition. Item 4 is a consequence of Item 3. Item 4
was the main idea behind the definition of $\mf nm$.
\end{proof}
The following proposition shows that $\mf{}{}$ is in fact a
differential algebra.
\begin{prop}
\label{16apr05}
For  $f\in \mf nm$ we have $\diff{f}\in \mf {n+1}{m+2}$ and
\begin{equation}
\label{khordim} \diff{(f||_mA)}= \diff{f}||_{m+2}A,\ \forall  A\in\GL(2,\R)
\end{equation}
\end{prop}
\begin{proof}
For $A\in \Mat_p(2,\Z)$ the term $\diff {(f||_mA)}$ is equal to:
\begin{eqnarray*}
 &= &
p^{m-n}\left (\sum_{i=0}^n \bn ni c_{A^{-1}}^i ((m-i)c_{A^{-1}}\ja(A,z)^{i-1-m}
f_i(Az)+\ja(A,z)^{i-m-2}\diff{f_i}(Az))\right ) \\
 &=&
p^{m-n}\left (\sum_{i=1}^{n+1}\bn {n}{i-1}c_{A^{-1}}^i\ja(A,z) ^{i-2-m}(m-i+1)f_{i-1}(Az) \right. \\ 
& + &  
\left. \sum_{i=0}^n \bn ni c_{A^{-1}}^i \ja(A,z)^{i-2-m}\diff{f_i}(Az)\right ) \\
&=&  p^{m-n}\left (\sum_{i=0}^{n+1}\bn {n+1}{i}c_{A^{-1}}^i\ja(A,z) ^{i-2-m}\tilde
f_i(Az)\right )
\end{eqnarray*}
where
$$\tilde f_i=\frac{i(m-i+1)}{n+1}f_{i-1}+
\frac{n+1-i}{n+1}\diff{f_i},\ i=0,1,\ldots,n+1, \ f_{-1}=f_{n+1}:=0
$$
Now $\diff{f}$ is a an $\mf{n+1}{m+2}$-function with the associated
$\mf{n-i}{m-2i}$ function $\tilde f_i$ for $i=0,1,2,\ldots,n+1$. The first item
of definition has been checked above, using the fact that $A\in\SL 2\Z$ and
$f||_mA=f$. We calculate $\tilde f_{ij}$ as above and then we get
\begin{eqnarray*}
\tilde f_{ij} &=&
\frac{i(m-i+1)}{n+1}f_{i-1+j}+\frac{n+1-i}{n+1}(\frac{j(m-2i-j+1)}{n-i+1}f_{i+j-1}+
\frac{n-i+1-j}{n-i+1}\diff{f_{i+j}})
\\  &= & \frac{(i+j)(m-(i+j)+1)}{n+1}f_{i+j-1}+
\frac{n+1-i-j}{n+1}\diff{f_{i+j}}=\tilde f_{i+j}
\end{eqnarray*}
This proves Item 2 of the definition of an $\mf{n+1}{m+2}$-function.
The third item can be checked using
$$
\diff{f}=2\pi i q\frac{df}{dq}
$$
\end{proof}
We consider $\mf{}{}$
as a graded algebra with $\deg(f)=m, \ f\in \mf nm$.
The nontrivial
statement about $\mf nm$-functions is Theorem \ref{23feb05} in
the Introduction. It can be interpreted also as the equality of
graded algebras $\mf{}{}=\C[g_1,g_2,g_3]$. It implies
that each $f\in\mf nm$ can be written as
$$
f=\sum_{i=0}^{n}f_ig_1^{i},\ f_i\in \mf 0{m-2i}
$$
We prove Theorem \ref{23feb05} in \S \ref{mainproof}.

The relations between the $g_i,i=1,2,3$ and their derivatives are given by
Ramanujan's equalities:
\begin{equation}
\label{ramanujan} \diff{g_1}=g_1^2- \frac{1}{12}g_2,\
\diff{g_2}=4g_1g_2- 6g_3,\ \diff{g_3}=6g_1g_3- \frac{1}{3}g_2^2
\end{equation}
(see for instance \cite{la95, nes01}). The proof of Theorem \ref{23feb05} will
contain a new proof of these equalities.
\subsection{Hecke operators}
\label{heop}
For $p\in\N$ let  $\SL 2\Z\backslash  \Mat_p (
2,\Z)=\{[A_1],[A_2],\ldots, [A_s]\}$. We define the $p$-th Hecke
operator in the following way
$$
T_pf=\sum_{k=1}^s f||_mA_k,\ \forall f\in \mf nm
$$
Lemma \ref{toulouse} implies that the above definition does not
depend on the choice of $A_k$ in the class $[A_k]$. Form Proposition
\ref{16apr05} one can deduce that the differential operator
$\diff{}$ commutes with the Hecke operator $T_p$.
\begin{prop}
\label{france}
 $T_p$ defines a map from $\mf nm$ to itself.
\end{prop}
This will be proved in \S \ref{heckeproof}.
One can take
$$
\tilde T_p:=\sum_{d\mid p, 0\leq b\leq d-1}\mat {\frac{p}{d}}b0d\in
\Z[\Mat_p(2,\Z)]
$$
and since for matrices $\mat ab0d$ the slash
operator $|_m$ is $p^n$ times $||_m$ we have
$T_pf=p^{-n}f|_m\tilde T_p$
and we get the expression (\ref{quedia}) in the Introduction.
In a similar way to the case of modular forms (see \cite{apo90} \S 6) one can
check that
$$
T_p\circ T_q=\sum_{d\mid (p,q)}d^{m-n-1}T_{\frac{pq}{d^2}}
$$
\subsection{The period domain}
The group $\SL 2\Z$ acts from left on the period domain $\pedo$ defined in
(\ref{perdomain}) and $G_0$ in (\ref{alggroup})
acts  from right. We consider holomorphic functions on
$$\L:=\SL 2\Z
\backslash \pedo
$$
as holomorphic functions
$$
f:\pedo\rightarrow \C,\hbox{ holomorphic satisfying },
f(Az)=f(z),\ \forall A\in\SL 2\Z, z\in\pedo
$$
The determinant function  is such a function.
The Poincar\'e upper half plane $\uhp$ is embedded in $\pedo$ in the following
way:
$$
z\rightarrow \tilde z=\mat{z}{-1}{1}{0}
$$
We denote by $\tilde\uhp$ the image of $\uhp$ under this map. For a
function $f$ on $\uhp$ we denote by $\tilde f$ the corresponding
function on $\tilde \uhp$.
\begin{prop}
There is a unique map  
$$
\phi: \mf{}{}\rightarrow {\cal O}(\pedo),\ f\mapsto \phi(f)=F
$$ of the
algebra of $\mf {}{}$-functions into the algebra of holomorphic functions on
$\pedo$ such that
\begin{enumerate}
\item
For all $f\in \mf {}{}$ the restriction of $F$ to $\tilde\uhp$ is
equal to $\tilde f$.
\item
For all $f\in \mf {}{}$ the holomorphic function $F$ is $\SL 2\Z$
invariant.
\item
We have
\begin{equation}
\label{7feb} F(x\cdot g)=k_2^nk_1^{n-m}\sum _{i=0}^n \bn ni
k_3^ik_2^{-i}F _i(x), \ \forall x\in\pedo,\ g\in G_0,
\end{equation}
where $F_i=\phi(f_i)$.
\end{enumerate}
Conversely, every holomorphic function $F$ on $\pedo$ which
satisfies 2,3  and whose restriction to $\tilde \uhp$ has finite
growth at infinity is of the form $F=\phi(f)$ for some $f\in \mf
nm$.
\end{prop}
We denote by $\check {\mf nm}$ the set of holomorphic functions on  $\pedo$
which restricted to $\tilde \uhp$ have finite growth at infinity,
are $\SL 2\Z$ invariant
and satisfy (\ref{7feb}).
For a classical modular form $f:\uhp\rightarrow \C$ of weight $m$ the
associated $F=\phi(f)$ is
$$
F(x)=x_3^m f(\frac{x_1}{x_3})\in \check{\mf 0m}
$$
We also have 
$$
\det \in \check{\mf 10}
$$
\begin{proof}
We have
\begin{equation}
\label{neeb}
\mat{x_1}{x_2}{x_3}{x_4}=\mat{\frac{x_1}{x_3}}{-1}{1}{0}
\mat{x_3}{x_4}{0}{\frac{\det(x)}{x_3}}
\end{equation}
So we expect $F$ to be defined by
$$
F(x)=F\left ( \mat{\frac{x_1}{x_3}}{-1}{1}{0}
\mat{x_3}{x_4}{0}{\frac{\det(x)}{x_3}}\right ):= x_3^{-m}\det(x)^n \sum
_{i=0}^n \bn ni x_4^ix_3^i\det(x)^{-i}f_i(\frac{x_1}{x_3})
$$
This $F$ restricted to $\tilde\uhp$ is $f$. By definition of $F$ one
can rewrite (\ref{4feb05}) in the form
$$
f(A\frac{x_1}{x_3})=(cx_1+dx_3)^{m-n}F\mat{x_1}{-d}{x_3}{c}
$$
where $A=\mat{a}{b}{c}{d}$.
We check item 3: Let
$$
g'=\mat {k_1'}{k_3'}{0}{k_2'}:=\mat{x_3}{x_4}{0}{\frac{\det(x)}{x_3}}
$$
\begin{eqnarray*}
\hbox{ RHS of (\ref{7feb}) } &= &
 k_2^nk_1^{n-m}\sum _{i=0}^n
\bn ni k_3^ik_2^{-i}F_i(x) \\
&= &
(k_2k_2')^n(k_1k_1')^{n-m}\sum _{i=0}^n\sum _{j=0}^{n-i}
\bn ni \bn {n-i}j k_3^ik_2^{-i}
k_2'^{-i}k_1'^{i}
k_3'^jk_2'^{-j}f_{ij}(\frac{x_1}{x_3}) \\
& =&
(k_2k_2')^n(k_1k_1')^{n-m}\sum _{r=0}^n\sum _{s=0}^{r}
\bn ns \bn {n-s}{r-s} k_3^sk_2^{-s}
k_2'^{-s}k_1'^{s}
k_3'^{r-s}k_2'^{-r+s}f_{s,r-s}(\frac{x_1}{x_3}) \\
&=&
(k_2k_2')^n(k_1k_1')^{n-m}\sum _{r=0}^n\sum _{s=0}^{r}
\bn nr  (k_2k_3'+k_3k_1')^r(k_2k_2')^{-r}f_r(\frac{x_1}{x_3}) \\
& = &
F(\mat{z}{-1}{1}{0}g'g)=F(xg) \\
\end{eqnarray*}
For all the equalities above we have used the same reasoning as in
the proof of  Lemma \ref{toulouse}.
We check that $F(Ax)=F(x),\ A\in \SL 2\Z$: The term $F(Ax)$ is equal to 
\begin{eqnarray*}
 &= &
(cx_1+dx_3)^{-m}\det(x)^n
\sum _{i=0}^n \bn ni (cx_2+dx_4)^i(cx_1+dx_3)^i\det(x)^{-i}
f_i(A\frac{x_1}{x_3}) \\
& =&
(cx_1+dx_3)^{-m}\det(x)^n \\ & . &
\sum _{i=0}^n \bn ni (cx_2+dx_4)^i(cx_1+dx_3)^i\det(x)^{-i}
(cx_1+dx_3)^{m-2i-(n-i)}F_i\mat{x_1}{-d}{x_3}{c} \\
& =&
F\left (\mat{x_1}{-d}{x_3}{c}
\mat{1}{\frac{cx_2+dx_4}{cx_1+dx_3}}{0}{\frac{\det(x)}{cx_1+dx_3}}\right )=F(x)
\end{eqnarray*}
Now let $F$ satisfy 2,3 and its restriction to $\tilde \uhp$ has
finite growth at infinity put $f=F\mid_{\tilde \uhp}$. First we note
that
$$
\mat abcd
\mat{z}{-1}{1}{0}=\mat{Az}{-1}{1}{0}\mat{\ja(A,z)}{-c}{0}{\ja(A,z)^{-1}\det(A)},\
A\in\GL (2,\R)
$$
Now
\begin{eqnarray*}
f(Az) &= & F\mat{Az}{-1}{1}{0} =F(\mat abcd
\mat{z}{-1}{1}{0}\mat{\det(A)\ja(A,z)^{-1}}{c}{0}{\ja(A,z)}) \\
 &=& \ja(A,z)^{n}\ja(A,z)^{m-n}\sum _{i=0}^n \bn ni
 c^i\ja(A,z)^{-i}f_i(x)=\sum _{i=0}^n \bn ni
 c^i\ja(A,z)^{m-i}f_i(x)
\end{eqnarray*}
This finishes the proof of our proposition.
\end{proof}
\subsection{Proof of Proposition \ref{france}}
\label{heckeproof}
We define
\begin{eqnarray*}
T_p:\check {\mf
nm}\rightarrow \check{\mf nm}, \ 
\check T_pF(x)=p^{m-2n-1}
\sum_{k=1}^s F(A_ix).
\end{eqnarray*}
This function has trivially its image in $\check{\mf nm}$. We
calculate the corresponding function in $\mf nm$: The term $T_pf(z)$ equals
\begin{eqnarray*}
 &= & p^{m-2n-1} \sum_{k=1}^s F(A_k\mat{z}{-1}{1}{0})
\\
& = &
p^{m-2n-1}\sum_{k=1}^s
 F\left (\mat{A_kz}{-1}{1}{0}\mat{\ja(A_k,z)}{-c}{0}{p.\ja(A_k^{-1},A_kz)}
\right ) \\
&= &
p^{m-2n-1}\sum_{k=1}^s (p.\ja(A_k^{-1},A_kz))^n(\ja(A_k,z))^{n-m}\sum
_{i=0}^n \bn ni
(-cp^{-1})^i\ja(A_k^{-1},A_kz)^{-i}f_i(A_kz) \\
&=&\sum_{k=1}^s  f||_m A_k
\end{eqnarray*}
where $A_k=\mat{a}{b}{c}{d}$. This proves Proposition \ref{france}.
\subsection{Some non-holomorphic functions on the period domain}
\label{chemishod}
On the complex manifold $\pedo$ we have the following $\SL 2\Z$ invariant
analytic functions:
$$
B_1:=\Im (x_1\overline{x_3}),\ B_2:=\Im (x_2\overline{x_4}),\
B_3:=x_1\overline{x_4}-x_2\overline{x_3}
$$
They define analytic functions on $\L$ which we denote them by the same letter.
They satisfy
$$
B_1\mid_{\tilde \uhp}(z)=\Im(z),\ B_1(xg)=B_1(x)|k_1|^2
$$
$$
B_2\mid_{\tilde \uhp}(z)=0,\ B_2(xg)=B_1(x)|k_3|^2+B_2(x)|k_2|^2+
\Im(B_3(x)k_3\overline{k_2})
$$
$$
B_3\mid_{\tilde \uhp}(z)=1,\ B_3(xg)=B_3(x)k_1\overline{k_2}+2\sqrt{-1}k_1
\overline{k_3}B_1(x).
$$
By equality (\ref{neeb}) one can easily see that
every point in $\pedo$ can be mapped to a point of
$\tilde \uhp$ by an action of a unique element of $G_0$.
This implies that  $\SL 2\Z$ invariant $B_i,\ i=1,2,3$, 
with the above properties are unique.
\subsection{Other topics}
\label{oldnew}
The literature of modular forms and its applications in number theory
is huge. The first question which naturally arises at this point is as
follows: Which part of the theory of modular forms can be generalized to
the context of differential modular forms and which arithmetic properties
can one expect to find? Since I am not expert in this area, I just mention
some subjects  which could fit well into this section.

The first of these, is the Eichler-Manin-Shimuara theory  of periods for cusp
forms (see \cite{HMM} and its references).
Note that the notion ``period'' in this theory, as far as I know,
has nothing to do with the notion of a period in this article. The notion
of period appears there because classical modular forms can be intrepreted
as sections of a tensor product of the cotangent bundle of a moduli curve and
hence a differential multi form, which can be integrated over some path in
the moduli curve (see \cite{sho80}). 
The differential modular forms are no longer interpreted
as sections of line bundles and this makes the situation more difficult.
Lewis type equations attached to differential modular forms will  be also 
of intereset (see \cite{HMM}).

Another theory which could be developed for differential modular forms is
Atkin-Lehner theory of old and new modular forms (see \cite{atle}).
This seems to me to be a quite accessible theory. The $L$-functions attached
to differential modular forms througth their Fourier expansion may also be of
interest (see \cite{bum97}).
\section{Families of elliptic curves and the Gauss-Manin connection}
\label{singular}
First, let us fix some notation.
For a ring $R$ we denote by $R[t]$ the polynomial
ring with coefficients in $R$ and the variable
$t:=(t_1,t_2,\ldots,t_s)$. We also define the affine space
$\A_R^s=\spec(R[t])$ defined over $R$. 
For simplicity we write $\A^s=\A_\C^s$. The set
of relative differential $i$-forms in $\A_R^s$ is denoted by 
$\Omega_{\A_R^s/\spec(R)}^i$.
\subsection{The family}
We consider the following family of elliptic curves
\begin{equation}
\label{family}
 {\cal E}_t: f=0
\end{equation}
$$
f:=y^2-4T_0(t)x^3+T_1(t)x^2+T_2(t)x+T_3(t)+T_4(t)y+T_5(t)xy,  \
t=(t_1,t_2,\ldots,t_s)\in T
$$
where
$$
T= \A^s\backslash \{\Delta=0\},\ T_i\in\C[t], \ i=0,1,\ldots,5
$$
and
{\tiny
$$
\Delta={\Big (}6912T_3^2-3456T_3T_4^2+432T_4^4{\Big )}T_0^4+
{\Big (}1152T_1T_2T_3-288T_1T_2T_4^2-576T_1T_3T_4T_5+
$$
$$
144T_1T_4^3T_5-256T_2^3+384T_2^2T_4T_5-288T_2T_3T_5^2-120T_2T_4^2T_5^2+
144T_3T_4T_5^3-4T_4^3T_5^3{\Big )}T_0^3+
(64T_1^3T_3-16T_1^3T_4^2-16T_1^2T_2^2+16T_1^2T_2T_4T_5-
$$
$$
48T_1^2T_3T_5^2+8T_1^2T_4^2T_5^2+8T_1T_2^2T_5^2-8T_1T_2T_4T_5^3+12T_1T_3T_5^4-
T_1T_4^2T_5^4-T_2^2T_5^4+T_2T_4T_5^5-T_3T_5^6 {\Big )}T_0^2
$$
}
The polynomial $\Delta\in\C[t]$ is the locus of parameters in which ${\cal E}_t$ is singular.
The reader is referred to \cite{mo} for the algorithms which
calculate $\Delta$ (see also the proof of Proposition \ref{8apr05}).
To make our notation simpler, we define
$$
\ring :=\C[t,\frac{1}{T_0}],\ \as_0:=\spec(R),\
\as_1:=\A^2\times (\A^s\backslash \{T_0=0\})=\spec(\C[x,y,t,\frac{1}{T_0}])
$$

\subsection{De Rham cohomology}
Let $\pi:\as_1\rightarrow \as_0$ be the
projection on the last $s$-coordinates. The following quotient
$$
H=\frac{\Omega^1_{\as_1}} { f\Omega^1_{\as_1}+df\wedge
\Omega^0_{\as_1}+ \pi^{-1}\Omega^1_{\as_0}\wedge
\Omega^0_{\as_1}+ d\Omega^0_{\as_1} } 
\cong
\frac{\Omega^1_{\as_1/\as_0}} {
f\Omega^1_{\as_1/\as_0}+df\wedge \Omega^0_{ \as_1/\as_0}+
d\Omega^0_{\as_1/\as_0} }
$$
is an $\ring$-module and will play the role of de Rham cohomology for us.
One may call $H$ the Brieskorn module associated to the family ${\cal
E}_t$ in analogy to the local modules introduced by Brieskorn in 1970.
The restriction of $H$ to the elliptic curve ${\cal E}_t$ gives us
$H_\dR^1({\cal E}_t)$. Set
\begin{equation}
\label{3mar05} \omega_1:=\frac{-2}{5}(2xdy-3ydx),\hbox{ and }
\omega_2:=\frac{-2}{7}x(2xdy-3ydx).
\end{equation}
\begin{prop}
The $\ring$-module $H$ is free and $\omega_1$ and
$\omega_2$ form a basis of $H$.
\end{prop}
\begin{proof}
We consider the classical Brieskorn module
$$
\tilde H=\frac{\Omega^1_{\as_1/\as_0}} { df\wedge
\Omega^0_{\as_1/\as_0}+ d\Omega^0_{\as_1/\as_0 } }
$$
It is a $\ring[f]$-module.
It is proved in \cite{mo} \S 6 Proposition 1 
that $\tilde H$ is freely
generated by $\omega_1,\omega_2$ as $\ring[f]$-module. Considering
the canonical map $\tilde H\rightarrow H$ we obtain the assertion of
our proposition.
\end{proof}
\subsection{Gauss-Manin connection}
Each element of the $\ring$-module $H$ can be interpreted as
a global section of the first cohomology bundle of the family ${\cal E}_t$
over $T$.
Since the Gauss-Manin connection is a connection
in the cohomology bundle, it is natural to find an algebraic definition for
it using the $\ring$- module structure of $H$. In this section we do this.
Define
$$
V:=\frac {\Omega^2_{\as_1}}
{df\wedge\Omega^1_{\as_1}+f\Omega^2_{\as_1}+
\pi^{-1}\Omega^1_{\as_0}\wedge \Omega^1_{\as_1}}
\cong \frac
{\Omega^2_{\as_1/\as_0}}
{df\wedge\Omega^1_{ \as_1/\as_0}+f\Omega^2_{ \as_1/\as_0}}
$$
\begin{prop}
\label{8apr05}
The polynomial $\Delta$ is a zero divisor of the $\ring$-module $V$,
i.e.
$$
\Delta\cdot V=0
$$
\end{prop}
\begin{proof}
We define
$$
\tilde V:=\frac {\Omega^2_{ \as_1/\as_0}}
{df\wedge\Omega^1_{ \as_1/\as_0}}
$$
and consider it as an $\ring[f]$-module. It is proved in \cite{mo} Lemma 4 
that
$B=\{dx\wedge dy,xdx\wedge dy\}$ form a basis of  $\tilde V$.
Let $A$ be the matrix of
multiplication by $f$ $\ring$-linear map in $\tilde V$ in the basis
$B$. Then $S(t,f):=\det(A-f.I_2)$ has the property
$S(t,f)\tilde V=0$, where $I_2$ is the identity two by two matrix.
This implies that $S(t,0)V=0$. Now $S(0,t)=6912.\Delta$ (In fact this is the 
way we have calculated $\Delta$).
\end{proof}
We have a well-defined differential map
$$
d: H\rightarrow V
$$
and we define the Gauss-Manin connection $H$ as follows:
$$
\nabla: H\rightarrow \Omega^1_{T}\otimes_{\ring} H
$$
$$
\nabla\omega=\frac{1}{\Delta}\sum_{i}\alpha_i\otimes\beta_i,\
\hbox{ where }\Delta d\omega=\sum_{i}\alpha_i\wedge \beta_i,\
\alpha_i\in \Omega^1_{\A^s},\ \beta_i\in\Omega^1_{\as_1}
$$
Let $U$ be an small open set in $U$ and $\{\delta_t\}_{t\in U},
\delta_t\in H_1({\cal E}_t,\Z)$ be a continuous family of
topological one dimensional cycles. The main property of the Gauss-Manin
connection is
\begin{equation}
\label{khodayakomak} d(\int_{\delta_t}\omega)=\sum \alpha_i
\int_{\delta_t}\beta_i,\ \nabla\omega=\sum_{i}\alpha_i\otimes
\beta_i,\ \alpha_i\in \Omega^1_{T},\ \beta_i\in H.
\end{equation}
 Let $\omega_1,\omega_2$ be a basis  of $H$ and define
$\omega=(\omega_1, \omega_2)^t$. The Gauss-Manin connection in this
basis can be written in the following way:
\begin{equation}
\label{violette}
\nabla\omega= A\otimes\omega,\
A=\frac{1}{\Delta}(\sum_{i=1}^s A_i dt_i)\in\Mat(2,\Omega^1_T),\
A_i\in\Mat (2,\C[t]).
\end{equation}
\begin{prop}
Let  $\tilde \omega=S\omega$ be another basis of $H$ and
$\nabla \omega=A\otimes\omega$.
 Then
$$
\nabla(\tilde\omega)=S(S^{-1}dS+A)S^{-1}\otimes \tilde\omega
$$
\end{prop}
\begin{proof}
We have
\begin{eqnarray*}
\nabla(\tilde\omega) &= &
\nabla(S\omega)=dS\otimes \omega+S\nabla\omega=dS. S^{-1}\otimes
\tilde\omega+SAS^{-1}\otimes
\tilde\omega \\
 &=& (dS.S^{-1}+SAS^{-1})\otimes\tilde\omega
\end{eqnarray*}
\end{proof}
\subsection{Classical differential forms}
\label{violtoul}
 Recall the canonical basis (\ref{3mar05}) of $H$.
Assume that in the family (\ref{family}) $T_3$ is of the form
$\tilde T_3+t_3$, where $\tilde T_3$ and all other $T_i$'s do not
depend  on $t_3$. Let $A$ be the matrix of the Gauss-Manin connection in
the basis $\omega$ as in (\ref{violette}). We are particularly
interested in $A_3$ in this case; it satisfies
\begin{equation}
\label{aemmm} \det(A_3)=321052999680\cdot T_0^4\cdot\Delta
\end{equation}
Define
$$\tilde \omega=\frac{1}{\Delta}A_3\omega$$ The classical way of
calculating Gauss-Manin connection leads to the equalities
\begin{equation}
\label{classical}
\tilde \omega_1=\frac{d\omega_1}{df}=\frac{dx}{y},\ \tilde \omega_2
= \frac{d\omega_2}{df}=\frac{xdx}{y}
\end{equation}
restricted to the elliptic curves ${\cal E}_t$ and up to exact
differential forms. Equality (\ref{aemmm}) implies that $\tilde
\omega$ form a basis of $H_\Delta$, the localization of $H$ over $\{1,
\Delta,\Delta^2,\ldots\}$.
\subsection{An example}
\label{examp}
The following example plays a basic role in the proof of Theorem \ref{23feb05}:
\begin{equation}
\label{4mar05}
 {\cal E}_t: \ f:=y^2-4t_0(x-t_1)^3+t_2(x-t_1)+t_3=0
\end{equation}
In this example
$$
\Delta=t_0(27t_0t_3^2-t_2^3)
$$
The calculation of the Gauss-Manin connection with respect to the canonical
basis (\ref{3mar05}) leads to:
{\tiny
$$
A_0=\mat {21/2t_0t_1t_2t_3-9t_0t_3^2+3/4t_2^3} {-21/2t_0t_2t_3}
{21/2t_0t_1^2t_2t_3+9t_0t_1t_3^2-1/2t_1t_2^3-5/8t_2^2t_3}
{-21/2t_0t_1t_2t_3-18t_0t_3^2+5/4t_2^3}
$$
$$
A_ 1 =\mat 0 0 {27t_0^2t_3^2-t_0t_2^3} 0
$$
$$
A_ 2 =\mat {-63/2t_0^2t_1t_3-5/4t_0t_2^2} {63/2t_0^2t_3}
{-63/2t_0^2t_1^2t_3+1/2t_0t_1t_2^2+15/8t_0t_2t_3}
{63/2t_0^2t_1t_3-7/4t_0t_2^2}
$$
$$
A_ 3 =\mat {21t_0^2t_1t_2+45/2t_0^2t_3} {-21t_0^2t_2}
{21t_0^2t_1^2t_2-9t_0^2t_1t_3-5/4t_0t_2^2}
{-21t_0^2t_1t_2+63/2t_0^2t_3}
$$
}
Now the same with respect to the classical basis $\tilde\omega$ of $H_\Delta$
is given by:
{\tiny
\begin{equation}
\label{rosa} A_0 =\mat {3/2t_0t_1t_2t_3-9t_0t_3^2+1/4t_2^3}
{-3/2t_0t_2t_3}
{3/2t_0t_1^2t_2t_3+9t_0t_1t_3^2-1/2t_1t_2^3+1/8t_2^2t_3}
{-3/2t_0t_1t_2t_3-18t_0t_3^2+3/4t_2^3}
\end{equation}
$$
A_1=\mat 0 0 {27t_0^2t_3^2-t_0t_2^3} 0
$$
$$
A_2 =\mat {-9/2t_0^2t_1t_3+1/4t_0t_2^2} {9/2t_0^2t_3}
{-9/2t_0^2t_1^2t_3+1/2t_0t_1t_2^2-3/8t_0t_2t_3}
{9/2t_0^2t_1t_3-1/4t_0t_2^2}
$$
$$
A_3 =\mat {3t_0^2t_1t_2-9/2t_0^2t_3} {-3t_0^2t_2}
{3t_0^2t_1^2t_2-9t_0^2t_1t_3+1/4t_0t_2^2}
{-3t_0^2t_1t_2+9/2t_0^2t_3}
$$
}
See \cite{mo, hos005} for the procedures which calculate all matrices above.
In the library {\tt brho.lib} the command {\tt gaussmaninp} calculate the
above matrix. Note that the canonical basis of the Brieskorn module 
in this library is $\mat{0}{\frac{-12}{5}}{\frac{-12}{7}}{0}\omega$.
\section{The period map}
\label{permap}
\subsection{Derivation of the period map}
Let $\omega=(\omega_1, \omega_2)^\tr$ be a basis of $H_\Delta$.
The period
map associated to the basis $\omega$ is given by:
$$
\pm: T\rightarrow \SL 2\Z\backslash \GL (2,\C),\ t\mapsto
\left [\frac{1}{\sqrt{2\pi i}}\mat
{\int_{\delta_1}\omega_1}
{\int_{\delta_1}\omega_2}
{\int_{\delta_2}\omega_1}
{\int_{\delta_2}\omega_2} \right ]
$$
It is well-defined and holomorphic. Here $\sqrt{i}=e^{\frac{2\pi
i}{4}}$ and $(\delta_1,\delta_2)$ is a basis of the $\Z$-module 
$H_1({\cal E}_t,\Z)$
such that the intersection matrix in this basis is
$\mat{0}{1}{-1}{0}$.

Let $\tilde \omega=S\omega, S\in
\Mat(2,\C[t,\frac{1}{T_0},\frac{1}{\Delta}])$ be another basis of
$H$ and $\tilde{\pm}$ be the associated period map. Then it is easy
to see that
\begin{equation}
\label{1apr05}
\tilde\pm=\pm \cdot S^\tr
\end{equation}
\begin{prop}
Let $\omega$ be a basis of $H_\Delta$ with
\begin{equation}
\label{29jan05} \nabla \omega=A \otimes
\omega,\  A\in \Mat (2,
\Omega^1_{T})
\end{equation}
Then
\begin{equation}
\label{4mar}
 d(\pm)(t)=\pm(t) \cdot A ^\tr,\ t\in T,
\end{equation}
where $d$ is the differential map.
\end{prop}
\begin{proof}
Let $\nabla_i:=\nabla_{\frac{\partial}{\partial t_i}}$ denote the Gauss-Manin connection with respect to
the parameter $t_i$. Then $\nabla_i \omega=A_i\omega$ and according
to (\ref{khodayakomak}) and (\ref{1apr05})
$$
\frac{\partial }{\partial t_i}\pm(t)=
\mat
{\int_{\delta_1}\nabla_i\tilde\omega_1}
{\int_{\delta_1}\nabla_i\tilde\omega_2}
{\int_{\delta_2}\nabla_i\tilde\omega_1}
{\int_{\delta_2}\nabla_i\tilde\omega_2}=\pm(t)\cdot A_i^\tr
$$
This proves the desired statement.
\end{proof}
For the basis $\omega=(\omega_1,\omega_2)^\tr$
if $\omega_1$ is such that its restriction to each elliptic curve
${\cal E}_t,\ t\in T$ is of the first kind then the period map is
defined  from $T$ into $\L:=\SL 2\Z\backslash \pedo$. For instance, the
period map associated to the classical basis (\ref{classical})
of $H_\Delta$ has this property.
\subsection{The Action of an algebraic group}
\label{mainproof}
We consider the family of elliptic curves (\ref{4mar05}).
It can be checked easily that (\ref{action}) is an action of $G_0$ on $\A^4$.
(This can be also  verified from the proof of the proposition bellow).
It is also easy to verify that $\A^4/G_0$ is isomorphic to $\Pn 1$ through
the map
\begin{equation}
s: \A^4/G_0\rightarrow \Pn 1,\  t\rightarrow [t_2^3:27t_0t_3^2-t_2^3]
\end{equation}
and so
\begin{equation}
\label{jinv}
j(t):=\frac{t_2^3}{27t_0t_3^2-t_2^3}
\end{equation}
is $G_0$-invariant and gives an isomorphy between $T/G_0$ and $\A$.
Recall the basis $\tilde \omega$ of $H_\Delta$ in (\ref{classical}).
\begin{prop}
\label{cano} The period $\pm$ associated to the basis $\tilde\omega$  is a
biholomorphism and
\begin{equation}
\label{gavril}
\pm(t\bullet g)=\pm(t)\cdot g,\ t\in\A^4,\ g\in G_0
\end{equation}
\end{prop}
\begin{proof}
We first prove (\ref{gavril}). Let
$$
\alpha: \A^2\rightarrow \A^2,\ (x,y)\mapsto
(k_2^{-1}k_1x-k_3k_2^{-1}, k_2^{-1}k_1^{2}y)
$$
Then
$$
k_2^{2}k_1^{-4}\alpha^{-1}(f)=y^2-4t_0k_2^{2}k_1^{-4} (
k_2^{-1}k_1x-k_3k_2^{-1}-t_1)^3+ t_2k_2^{2}k_1^{-4}(
k_2^{-1}k_1x-k_3k_2^{-1}-t_1)+t_3k_2^{2}k_1^{-4}
$$
$$
y^2-4t_0k_1^{-1}k_2^{-1}(x-(t_1k_2k_1^{-1}+k_3k_1^{-1}))^3+
t_2k_1^{-3}k_2(x-(t_1k_2k_1^{-1}+k_3k_1^{-1}))+t_3k_1^{-4}k_2^{2}
$$
This implies that $\alpha$ induces an isomorphism of elliptic curves
$$
\alpha: {\cal E}_{t\bullet g}\rightarrow {\cal E}_t
$$
Now
$$
\alpha^{-1} \tilde\omega= \mat{k_1^{-1}}{0}{-k_3k_2^{-1}k_1^{-1}}{k_2^{-1}}\tilde\omega=\mat {k_1}{0}{k_3}{k_2}^{-1} \tilde\omega
$$
By the equality (\ref{1apr05}) we have
$$
\pm(t)= \pm(t\bullet g).g^{-1}
$$
which proves (\ref{gavril}).

 The matrix of the Gauss-Manin connection
in the basis $\tilde\omega$ for the family (\ref{4mar05}) is calculated
in \S \ref{examp}. Let $B$ be a $4\times 4$
matrix and the $i$-th row of $B$ constitutes of the first and
second rows of $A_i$. Then
$$
\det(B)=\frac{3}{4}t_0\Delta^3
$$
shows that the period map $\pm$ is regular at each point $t\in
T$ and hence it is locally a biholomorphism.


The period map $\pm$ induces a local biholomorphic map
$\bar \pm: T/G_0\rightarrow \SL 2\Z\backslash \uhp\cong \C$ and so we have the
local biholomorphism $\bar \pm \circ j^{-1} : \A\rightarrow \A$.
One can compactify $\SL 2\Z\backslash \uhp$ by adding the cusp $\SL
2\Z/\Q=\{c\}$ (see \cite{miy}) and the map $ \bar \pm \circ j^{-1}$ is
continuous at $v$ sending $v$ to $c$, where $v$ is the point induced by
$t_027t_3^2-t_2^3=0$ in $\A^4/G_0$.
 Using Picard's Great Theorem we conclude that $j^{-1}\circ \bar \pm$ is a
biholomorphism and so $\pm$ is a biholomorphism.
\end{proof}
\subsection{The inverse of the period map}
We denote by
$$
F=(F_0,F_1,F_2,F_3):\pedo \rightarrow T
$$
the inverse of the period map.
\begin{prop}
The following is true:
\begin{enumerate}
\item
$F_0(x)=\det(x)^{-1}$.
\item
For $i=2,3$
$$
F_i=\det(x)^{1-i}  \check{\es{ i}}\in \check{\mf 0{2i}}
$$
where $\es {i}$ is the Eisenstein series (\ref{eisenstein}).
\item
$F_1=\check{g_1} \in \check {\mf {1}{2}}$.
\end{enumerate}
\end{prop}
\begin{proof}
Taking $F$ of (\ref{gavril}) we have
$$
F_0(xg)=F_0(x)k_1^{-1}k_2^{-1}, 
$$
\begin{equation}
\label{2apr05}
F_1(xg)=F_1(x)k_1^{-1}k_2+k_3k_1^{-1},\
\end{equation}
$$
F_2(xg)=F_2(x)k_1^{-3}k_2,\  F_3(xg)=F_3(x)k_1^{-4}k_2^2,\ \forall 
x\in\L,\
g\in G_0
$$
By Legendre's  Theorem $\det(x)$ is equal to one on $V:=\pm 
(1\times 0\times \A\times \A)$ and so the same is true for
$F_0\det(x)$. But the last function is invariant under the action of
$G_0$ and so it is the constant function $1$.  This proves the first
item.
Let $G_i=F_i\det(x)^{i-1},\ i=1,2,3$. The equalities 
(\ref{2apr05}) imply that $G_i, i=2,3$ do not depend on $x_2,x_4$.
Now the map $(t_2,t_3)\rightarrow \pi\circ \pm(1,0,t_2,t_3)$, where $\pi$ is 
the projection on the $x_1,x_3$ coordinates, is the classical period map
(see for instance see  \cite{sai01} and  and the appendix of 
\cite{ka73}) and this implies that 
$G_i= \check {\es{i}},\ i=2,3$. 
Note that in our definition of the 
period map the factor $\frac{1}{\sqrt{2\pi i}}$ appears.
In particular $F_i, i=2,3$  have finite growth at infinity. 
The fact that $F_1$ has finite growth at infinity follows form
the Ramanujan relations (\ref{ramanujan}) and the equality 
$\frac{d}{dz}=2\pi i q \frac{d}{dq}$.
Since $G_1\in \check {\mf 12}$, $\check {\mf 12}$ is
a  one dimensional
space, both $g_1,G_1$ satisfy (\ref{2apr05})
and ${\mf 0{2}}=\emptyset$, we have $G_1=g_1$. 
\end{proof}
\subsection{Ramanujan  relations}
We proved in Lemma \ref{cano} that the period map $\pm$ associated
to $\tilde \omega$ is a biholomorphism. According to (\ref{4mar}), the
inverse $F$ of $\pm$ satisfies the differential equation
$$
x.A(F(x))^\tr=I
$$
We consider $\pm$ as a map sending the vector $(t_0,t_1,t_2,t_3)$ to
$(x_1,x_2,x_3,x_4)$. Its derivative at $t$ is a $4\times 4$ matrix whose
$i$-th column constitutes of the first and second row of
$\frac{1}{\Delta}xA_i^\tr$. We  use (\ref{rosa})
to derive the equality {\tiny
$$
(dF)_x=(d\pm)_t^{-1}=
$$
$$
\det(x)^{-1}\left ( \begin{array}{cccc} -F_0x_4 &F_0x_3
&F_0x_2 &-F_0x_1
\\ \frac{1}{12F_0}(12F_0F_1^2x_3-12F_0F_1x_4-F_2x_3)
&-F_1x_3+x_4 &\frac{1}{12F_0}(-12F_0F_1^2x_1+12F_0F_1x_2+F_2x_1)
&F_1x_1-x_2
\\4F_1F_2x_3-3F_2x_4-6F_3x_3
&-F_2x_3 &-4F_1F_2x_1+3F_2x_2+6F_3x_1 &F_2x_1
\\\frac{1}{3F_0}(18F_0F_1F_3x_3-12F_0F_3x_4-F_2^2x_3)
&-2F_3x_3 &\frac{1}{3F_0}(-18F_0F_1F_3x_1+12F_0F_3x_2+F_2^2x_1)
&2F_3x_1
\end{array} \right ).
$$
} 
For $g_i:=F_i\mid_{\tilde \uhp}$ we obtain the Ramanujan
relations (\ref{ramanujan}).
\subsection{Other topics}
 We have seen that the period map satisfies the
differential equation
\begin{equation}
\label{pfeq}
d(\pm^\tr)=A\pm^\tr
\end{equation}
This can be interpreted as a multi-variable Fuchsian system/Picard-Fuchs
equation.
Restricting (\ref{pfeq}) to a line in $ \A^4$, one gets a usual
one variable Picard-Fuchs equation. For many of them the solution and hence
the associated abelian integral can be given explicitly by
hypergeometric functions (see \cite{beu95}). It seems to me that with
the notion of period map in this article it is possible to classify all the
algebraic values of $G$-functions obtained by the methods in \cite{beu95}.

Another interesting subject could be the construction of
a logarithmic  structure on $T$ in the context of algebraic geometry.
Such constructions for the moduli of polarized Hodge structures is done in
\cite{kaus00} and \cite{kaus04} in the context of analytic geometry .

\section{The Ramanujan foliation}
\label{ramanfoli}
The whole theory developed in the previous sections could be done by
setting $t_0=1$. In this section we assume that $t_0=1$ and use the same
notations for $\pm,\pedo,G_0,T, \Delta$ and so on. For instance
redefine
$$
\pedo:=\{x=\mat {x_1}{x_2}{x_3}{x_4}\in \GL(2,\C)\mid \Im(x_1\ovl{x_3})>0,\
\det(x)=1 \}
$$
and
$$
G_0=\{\mat{k}{k'}{0}{k^{-1}}\mid \ k'\in\C, k\in \C^*\}
$$
The action of $G_0$ on $\A^3$ is given by
$$
t\bullet g:=(
 t_1k^{-2}+k'k^{-1},
t_2k^{-4}, t_3k^{-6}), t=(t_1,t_2,t_3)\in\A^3,
g=\mat {k}{k'}{0}{k^{-1}}\in G_0
$$
We also define
$$
g=(g_1,g_2,g_3): \uhp\rightarrow T\subset \A^3
$$
\subsection{The Ramanujan foliation}
We write (\ref{raman}) in the vector filed form
$$
\Ra:=
(t_1^2-\frac{1}{12}t_2)\frac{\partial }{\partial t_1}+
( 4t_1t_2-6t_3)\frac{\partial }{\partial t_2}+
(6t_1t_3 -\frac{1}{3}t_2^2)\frac{\partial }{\partial t_3}
$$
It is also useful to define the differential forms
$$
\eta_{1}:=(t_1^2-\frac{1}{12}t_2)dt_2-
( 4t_1t_2-6t_3)dt_1
$$
$$
\eta_2:=( 4t_1t_2-6t_3)dt_3-
(6t_1t_3 -\frac{1}{3}t_2^2)dt_2,
$$
$$
\eta_3:= (t_1^2-\frac{1}{12}t_2)
dt_3-(6t_1t_3 -\frac{1}{3}t_2^2)dt_1
$$
and say the the foliation $\F(\dR)$ is induced by $\eta_i,i=1,2,3$.
The singularities of (\ref{raman}) are the points $t\in\A^3$ such
that $\Ra(t)=0$. It turns out that
$$
\sing(\Ra)=\{(t_1,12t_1^2,8t_1^3)\mid t_1\in \A \},
$$
which is a one dimensional curve and lies in $\{\Delta=0\}$.
We have
\begin{eqnarray*}
d\Delta(\Ra) & = & (2.27t_3dt_3-3t_2^2dt_2)(\Ra) \\
 & =& 2.27t_3 (6t_1t_3-\frac{1}{3}t_2^2)-3t_2^2(4t_1t_2-6t_3) \\
 &= & 12t_1\Delta
\end{eqnarray*}
This implies that the variety
$\Delta_0:=\{\Delta=0\}$ is invariant by the foliation $\F(\Ra)$. Inside
$\Delta_0$ we have the algebraic leaf $\{(t_1,0,0)\in\A^3\}$ of $\F(\Ra)$.
On $\sing(\Ra)$ we have a special point $p_\infty$ given by (\ref{pinfty}).
It is the limit of $g(z)$ when $\Im(z)$ tends to $+\infty$.
\begin{prop}
\label{16apr}
The following is a uniformization of the foliation $\F(\Ra)$
restricted to $T$:
$$
u: \uhp\times (\A^2\backslash \{(0,0)\})\rightarrow T
$$
\begin{equation}
\label{9apr05}
(z,c_2,c_4)\rightarrow g(z)\bullet
\mat{(c_4z-c_2)^{-1}}{c_4}{0}{c_4z-c_2}=
\end{equation}
$$
(g_1(z)(c_4z-c_2)^2+(c_4z-c_2), g_2(z)(c_4z-c_2)^4, g_3(z)(c_4z-c_2)^6).
$$
\end{prop}
\begin{proof}
One may check directly that for fixed $c_2,c_4$ the map induced by $u$
is tangent to (\ref{raman}).
In general it is not a solution of (\ref{raman}). This is the main reason
for naming ``foliation''. We give another proof which uses the period map:
From (\ref{rosa}) we have
$$
d(\pm)(t)=\pm(t)\mat {\frac{3}{4}\eta_2}{\frac{3}{2}(3t_3dt_2-2t_2dt_3)}
{\frac{9}{2}t_3\eta_1-3t_2\eta_3+\frac{3}{2}t_1\eta_2}
{-\frac{3}{4}\eta_2 }^{\tr}.
$$
Therefore
$$
d(\pm(t))(\Ra(t))=\pm(t)\mat{0}{0}{*}{0}=\mat{*}{0}{*}{0}.
$$
This implies that the $x_2$ and $x_4$ coordinates of
the pull forward of the vector field $\Ra$ by $\pm$ are zero. Therefore, the
leaves of $\F(\Ra)$ in the period domain are of the form
$$
\mat{z(c_4z-c_2)^{-1}}{c_2}{(c_4z-c_2)^{-1}}{c_4}=
\mat{z}{-1}{1}{0}
\mat{(c_4z-c_2)^{-1}}{c_4}{0}{c_4z-c_2}.
$$
\end{proof}
Due to Proposition \ref{16apr}, 
our foliation $\F(\Ra)$ can be considered  as a kind of 
Hilbert modular foliation (see \cite{peme}). 
\subsection{The family $y^2-4x^3+t_1x^2+t_2x+t_3$}
The family (\ref{khodaya}) can be rewritten in the
form
$$
y^2-4t_0x^3+12t_0t_1x^2+(-12t_0t_1^2+t_2)x+(4t_0t_1^3-t_2t_1+t_3)=0
$$
The mapping $$ \alpha: \A^4\rightarrow \A^4, t\mapsto (t_0,
12t_0t_1, -12t_0t_1^2+t_2,  4t_0t_1^3-t_2t_1+t_3) $$ is an
isomorphism and so we can re state Proposition \ref{cano} for the family
$$
{\cal E}_t: y^2-4t_0x^3-t_1x^2-t_2x-t_3=0.
$$
The inverse of the period map in this case is given by
$G=(G_0,G_1,G_2,G_3)$ with
$$
G_0=F_0\in \check {\mf {-1}{0}},\  G_1=12F_0F_1\in \check {\mf
{0}{2}}, $$ $$
G_2=-12F_0F_1^2+F_2\in \check {\mf {1}{4}}  ,\
G_3=4F_0F_1^3-F_2F_1+F_3\in \check {\mf {2}{6}}
$$
In this case the singular fibers are parameterized by the zeros of
$$
\Delta:=
t_0
(432t_0^2t_3^2+72t_0t_1t_2t_3-16t_0t_2^3+4t_1^3t_3-t_1^2t_2^2).
$$
The Ramanujan relations take the simpler form:
\begin{equation}
\label{ramanu}
\left \{ \begin{array}{l}
\dot {t_1}= -t_2 \\
\dot{t_2}= -6t_3  \\
\dot {t_3}= t_1t_3-\frac{1}{4}t_2^2
\end{array} \right.
\end{equation}
They have the solution
$$
g:=(12g_1, -12g_1^2+g_2,4g_1^3-g_2g_1+g_3)
$$
\section{Proofs}
\label{proofs}
Now we are in a position to prove the theorems announced in the Introduction.
\subsection{Proof of Theorem \ref{23feb05}}
We prove that $\check{\mf {}{}}$ as a $\C$-algebra is freely generated by
$\frac{1}{F_0}, F_i,\ i=0,1,2,3$. 
Let $\tilde F\in\check{\mf{n}{m}}$ and 
$\tilde F_i\in \check{\mf {n-i}{m-2i}}$ be its 
associated functions. Since the period map $\pm: T\rightarrow \L$ is a 
biholomorphism, there exist holomorphic functions 
$p_i \ i=0,1,\ldots,n,\ p_0:=p$ defined on $T$ such that 
$\tilde F_i=p_i(F_0,F_1,F_2,F_3)$.
The property (\ref{7feb}) of $\tilde F$ implies that:
\begin{equation}
\label{jabr}
p(t\bullet g)=k_{2}^nk_1^{n-m}\sum_{i=0}^n\bn nik_3^ik_2^{-i}p_i(t),\ 
\forall g\in G_0,\ t\in T.
\end{equation}
Take $g=\mat{1}{t_1}{0}{1}$ and $t=(t_0,0,t_1,t_3)$. Then
$$
p(t_0,t_1,t_2,t_3)=\sum_{i=0}^n\bn nit_1^ip_i(t_0,0,t_2,t_3).
$$
This implies that $p$ is a polynomial of degree at most $n$ in 
the variable $t_1$. Let us re write $p(t)=\sum_{i=0}^{n} t_1^i q_i$, 
where $q_i$'s are holomorphic functions on $\A^3\backslash
\{t\in\A^3\mid \Delta=0\}$. We apply (\ref{jabr}) to 
$g=\mat{k}{0}{0}{t_0^{-1}k^{-1}}$ and consider the coefficients of 
$t_1^i,\ i=1,2\ldots,n$. We get 
\begin{equation}
\label{13apr}
q_i(1,t_2t_0^{-1}k^{-4}, t_3t_0^{-1}k^{-6})=t_0^{-n}k^{-m+2i}q_i(t),\ i=1,2,\ldots,n. 
\end{equation}
Take $t_0=1$. 
The growth condition on $\tilde F$ is translated through the period map 
into the following fact: $p$ restricted to 
a transversal disk to $\Delta=0$ at $p_\infty$ has  a holomorphic extension 
to $p_\infty$. This will also imply the similar growth conditions for
$q_i$'s. The classical fact that the set of modular forms is generated by
the Eisenstein series $\es{2}$ and $\es{3}$ and (\ref{13apr}) with $t_0=1$ 
imply that $q_i(1,t_2,t_3)$ 
is a homogeneous polynomial of degree $m-2i$ in the graded ring 
$\C[t_2,t_3],\ \deg(t_2)=4,\deg(t_3)=6$. 
We conclude that $p$ is of the form
$$
p=t_0^n \sum_{i=0}^nt_1^i q_i(1, t_2t_0^{-1}, t_3t_0^{-1}).
$$
In other words, $p$ is homogeneous of degree $m$ in the variables 
$t_1,t_2,t_3$ with 
$\deg(t_i)=2i,\ i=1,2,3$.
\subsection{Proof of Theorem \ref{realanal}}
In \S \ref{chemishod} we described some analytic functions $B_i,\ i=1,2,3$,
on $\L$ which  had some compatibility properties with the action of $G_0$ on
$\L$. We use Proposition \ref{cano} and transfer them to the world of
coefficients. This will prove the existence and uniqueness of 
the functions $B_1,B_2,B_3$.
For the sake of simplicity we have used the same letters to name these
functions.

By the properties that $B_1$ has we can say more about it.
In (\ref{B1}) we put $k=1$ and we conclude that $B_1$ is independent
of the variable $t_1$. The function
$
B_2\cdot |\Delta|^\frac{1}{6}
$
is $G_0$ invariant and so there is an analytic function
$b_2:\A\rightarrow \R$
such that
$$
B_2(t)=\frac{b_2(j(t))}{|\Delta(t)|^\frac{1}{6}}.
$$
Translating this to $\uhp$, we have
$$
\Im(z)=\frac{b_2(j(z))}{|\Delta(z)|^\frac{1}{6}}
$$
where the above $j$ and $\Delta$ are the ones on \S \ref{algebra}.

The proof of the last part of the theorem is as follows: On $M_0$ an
$x\in \pedo$ can be written in the form $\mat
{x_1}{x_4r}{x_3}{x_4},\ r\in\R, x_4(x_1-rx_3)=1$. Then
\begin{equation}
\label{10apr05}
B_3(x)=\overline{x_4}(x_1-rx_3)=\frac{\overline{x_4}}{x_4}.
\end{equation}
\subsection{Proof of Theorem \ref{foli}}
We follow the notations introduced in \S \ref{ramanfoli}. In
particular we work with the family (\ref{khodaya}) with $t_0=1$. The
leaves of the pull-forward of the foliation $\F(\Ra)$ by the period
map $\pm$ have constant $x_2$ and $x_4$ coordinates. By definition
of $B_2:=\Im(x_2\overline{x_4})$ in the period domain, we conclude that
$M_r$'s are $\F(\Ra)$-invariant. The equality (\ref{10apr05})
implies that $N_w$'s are $\F(\Ra)$-invariants.

Let us now prove item 2: Take $t\in K$ and a cycle $\delta\in
H_1({\cal E}_t,\Z)$ such that $\int_\delta\tilde\omega_2=0$ and
$\delta$ is not of the form $n\delta'$ for  some $2\leq n\in \N$ and
$\delta'\in H_1({\cal E}_t,\Z)$. We choose another $\delta'\in
H_1({\cal E}_t,\Z)$ such that $(\delta',\delta)$ is a basis of
$H_1({\cal E}_t,\Z)$ and the intersection matrix in this basis is
$\mat{0}{1}{-1}{0}$. Now $\pm(t)$ has zero $x_4$-coordinates and so
its $B_2$ is zero. This implies that $K\subset M_0$. It is dense because
an element $\mat{x_1}{x_4r}{x_3}{x_4}\in M_0\subset \L$ 
can be approximated by the elements in $M_0$ with $r\in \Q$.  

The image of the map $g$ is the locus of the points $t_0$ in $T$
such that $\pm(t_0)$ is of the form $\mat{z}{-1}{1}{0}$ in a basis
of $H_1({\cal E}_t,\Z)$. We look $g$ as a function of $q=e^{2\pi
iz}$ and we have
$$
g(0)=p_\infty,\ \frac{\partial g}{\partial
q}(0)=(-24a_1,240a_2,-504a_3)
$$
where $p_\infty=(a_1,a_2,a_3)$ is the one in (\ref{pinfty}). This
implies that the image of $g$ intersects $\sing(\Ra)$ transversally.
For a $t\in K$ the $x_4$-coordinate of $\pm$ is zero and the leaf
through $t$, namely $L_t$, has constant $x_2$-coordinate, namely
$c_2$. By (\ref{9apr05}) $L_t$ is uniformized by
$$
u(z)=(c_2^2g_1(z), c_2^4g_2(z),c_2^6g_3(z)),\ z\in\uhp .
$$
This implies that $L_t$ intersects $\sing(\Ra)$ transversally at
$(c_2^2a_1,c_2^4a_2,c_2^6a_3)$.

We prove item 3: Let $t\in T$ and the leaf $L_t$ through $t$ have an
accumulation point at $t_0\in T$. We use the period map $\pm$ and
look $\F(\Ra)$ in the period domain. For $(c_2,c_4)\in
\A^2\backslash \{0\}$ the set $S=\{A(c_2,c_4)^\tr\mid A\in \SL
2\Z\}$ has an accumulation point in $\A^2$ if and only if
$\frac{c_2}{c_4}\in\R\cup{\infty}$ or equivalently $B_2(t)=0$.
\subsection{Proof of Theorem  \ref{ellip}}
Let $k$ be an algebraically closed field of charachteristic $0$ ,
for instance take $k=\overline{\Q}$. 
\begin{prop}
\label{classspace}
The quasi affine variety
$$
T=\spec(k[t, \frac{1}{\Delta}])
$$
is the moduli of $(E, [\omega_1],[\omega_2])$'s, where $E$ is an elliptic
curve defined over $k$, $\omega_1$ is a differential form of the first kind 
on $E$ and $([\omega_1], [\omega_2])$ is a basis of $H^1_\dR(E)$.
\end{prop}
\begin{proof}
For simplicity we do not write more $[.]$ for differential forms.  
The $j$ invariant (\ref{jinv}) classifies the ellipric curves over $k$ (see 
\cite{har77} Theorem 4.1). Therefore, for a given elliptic curve $E/k$ we can
find parameter $t\in \A_k^4$ such that $E\cong {\cal E}_t$ over $k$. 
Under this isomorphism we write
$$
\matt{\omega_1}{\omega_2}=g^\tr\matt{\frac{dx}{y}}{\frac{xdx}{y}},\ 
\hbox{ in }\ H^1_\dR({\cal E}_t)
$$
for some $g\in G_0$, where $\omega_1,\omega_2$ are as in the proposition. 
Now, the triple $(E,\omega_1,\omega_2)$ is isomorphic to 
$({\cal E}_{t\bullet g}, \frac{dx}{y}, \frac{xdx}{y})$.
Since $j: \A^4/G_0\rightarrow \A$ is an isomorphism, every triple
$(E,\omega_1,\omega_2)$ is represented exactly by one parameter $t\in T$. 
\end{proof}
By Proposition \ref{classspace} the hypothesis of Theorem \ref{ellip} gives 
us a parameter $t\in T$ such
that $\int_{\delta}\frac{xdx}{y}=0$, for some $\delta\in H_1({\cal
E}_t,\Z)$. We can assume that $\delta$ is not a multiple of another cycle
in $H_1({\cal E}_t,\Z)$. 
The corresponding period matrix of $t$ in a basis
$(\delta',\delta)$ of $H_1({\cal E}_t,\Z)$ has zero $x_4$-coordinate
and so the numbers
$$
t_0=\det(x)^{-1},\
t_i=\det(x)^{1-i}x_3^{-2i}\es{i}(\frac{x_1}{x_3}),\ i=2,3,\
t_1=F_1(\mat{x_1}{x_2}{x_3}{0})=\det(x)x_3^{-2}\es{1}(\frac{x_1}{x_3})
$$
all are in $\overline{\Q}$. 
This implies that for $z=\frac{x_1}{x_3}\in\uhp$ we have
$$
\frac{\es{3}}{\es{1}^3}(z), \frac{\es{4}}{\es{1}^2}(z),\
\frac{\es{3}^2}{\es{2}^3}(z)\in \overline{\Q}.
$$
This is in contradiction with

{\bf Theorem} (Nesterenko 1996, \cite{nes01})
{\it For any $z\in\uhp$, the set
$$
e^{2\pi i z},\ \frac{g_1(z)}{a_1}, \frac{g_2(z)}{a_2},\frac{g_3(z)}{a_3}
$$
contains at least three algebraically independent numbers over $\Q$. }

A direct corollary of Theorem \ref{ellip} is that the multivalued
function
$$
I(t)=\frac{\int_{\delta_t}\frac{xdx}{y}}{
\int_{\delta_t}\frac{dx}{y} }
$$
defined in $T$ never takes algebraic values for algebraic $t$.
\subsection{Other topics}
As a person who has started his mathematical career by studying
holomorphic foliations in complex manifolds, I would be interested
to describe completely the dynamics of the foliation
$\F(\Ra)$ and in particular the behavior of the leaves near
$\Delta=0$. The leaves of $\F(\Ra)$ in $\Delta=0$ parameterize
degenerated elliptic curves. Can one describe their  behavior by
abelian integrals?

Our proof of Theorem \ref{23feb05} is completely based on the
existence of the nice family (\ref{khodaya}) for which the period map
is a biholomorphism. To prove Theorem \ref{23feb05} for certain
subgroups of $\SL 2\Z$ by the methods of this article, one must find
a four parameter family of elliptic curves such that the period map
is an etale covering, i.e. it is a  local biholomorphism of finite
degree. The inverse of the period map is a multi valued function
whose restriction to the simply connected domain $\tilde \uhp$ gives
rise to a finite number of holomorphic functions on $\uhp$. These
new functions  are differential modular with respect to a subgroup
of $\SL 2\Z$ which can be calculated by some topological data, such as
the homotopy group of $T$, attached to the period map. 

The Moduli space $T$ is not a Shimura variety (see \cite{mil03}). 
In this point it would be too much speculation to say that spaces like
$T$ can be constructed for arbitrary moduli of polarized Hodge structures.
Nevertheless, more constructions of such spaces may result in generalizations
of Shimura varieties.
{\tiny

\def\cprime{$'$} \def\cprime{$'$}

\bibliographystyle{plain}
}
\end{document}